\documentclass[a4paper]{article}
\usepackage{amsmath}
\usepackage[colorlinks=true,breaklinks=true,citecolor=red,linkcolor=blue,urlcolor=red,unicode]{hyperref}
\usepackage[hyperpageref]{backref}
\renewcommand*{\backref}[1]{}
\renewcommand*{\backrefalt}[4]{\ifcase #1\or(p.~#2).\else(pp.~#2).\fi}
\renewcommand*{\backreftwosep}[1]{ and }
\renewcommand*{\backreflastsep}[1]{ and }

\newlength{\marge}
\paperwidth
\usepackage{amssymb}
\usepackage{mathrsfs}
\usepackage{textgreek}

\newtheorem{thm}{Theorem}[section]
\newtheorem{thmf}[thm]{Th\'eor\`eme}
\newtheorem{lem}[thm]{Lemma}

\newtheorem{prp}[thm]{Proposition}

\newtheorem{dfn}[thm]{Definition}

\newtheorem{construction}[thm]{Construction}
\newtheorem{Q}[thm]{Question}

\def\dem{{\it Proof.\/ }} 
\def\rem{\refstepcounter{thm}
{\bf Remark \thethm\hskip1ex}}

\def\exa{
\addtocounter{thm}{1}{\bf Example \thethm}\hskip1ex}

\def\vskipa{\vskip\abovedisplayskip}
\def\vskipb{\vskip\belowdisplayskip}
\def\ie{{\it i.\,e.}\ } 
\def\vs{{\it vs.}\ }

\def\toll#1{\vtop{\ialign {##\crcr \rightarrowfill \crcr 
\noalign{\kern -1pt \nointerlineskip \vskip2pt
} 
$\hfil \scriptstyle{\ #1\ } \hfil $\crcr }}}

\def\({{\fam0\rm (}}
\def\){\/{\fam0\rm )}} 
\def\]{\mathopen]}
\def\[{\mathclose[}
\mathchardef\hook="312C

\def\tol{\mathop{-\mkern-4mu\smash-\mkern-4mu\smash\longrightarrow}\limits}

\def\imp{\Rightarrow}

\def\ecks{\rule{2mm}{2mm}}
\def\bull{\leavevmode\kern .1ex\vrule height 1ex width .9ex depth
-.1ex \kern .8ex} 
\def\bloc{\bull}
\def\eck{\nolinebreak\hspace{\fill}\ecks}
\def\barre{\rule[-2.5pt]{0.8pt}{10pt}}
\def\lmes{\kern0.7pt\barre\kern0.7pt}
\def\rmes{\kern0.7pt\/\barre\kern0.7pt}
\def\mes#1{\mbox{$\lmes #1 \rmes$}}

\def\bi#1#2{\bigg(\kern-0.5pt{{#1}\atop{#2}}\kern-0.5pt\bigg)}
\def\bip#1#2{\left(\kern-3pt {{#1}\atop{#2}}\kern-2pt \right)}
\def\se{\subseteq}
\def\HM#1{\setbox0=\hbox{#1}\dimen0=\wd0 #1 \kern-\dimen0
\setbox0\hbox{\raise3pt\hbox{$^\frown$}}\advance\dimen0 by -\wd0
\box0\kern\dimen0}

\def\P#1{\mathop{\mathbb P}\nolimits
\left[\kern0.5pt{#1}\kern0.5pt\right]}
\def\E{\mathop{\mathbb E}\nolimits} 
\def\e{\mathop{\fam0 e}\nolimits}

\def\ii{{\fam0 i}}

\def\Zeta{{\fam0 Z}} 
\def\eps{\varepsilon}

\def\Z{\hbox{${\mathbb Z}$}}

\def\T{\hbox{${\mathbb T}$}}
\def\N{\hbox{${\mathbb N}$}}

\def\C{\hbox{${\mathbb C}$}}

\def\ds{\displaystyle} 

\def\lst#1#2{#1_1,\dots,#1_{#2}}

\def\sm#1#2{#1_1+\dots+#1_{#2}} \def\bip#1#2{\left(\kern-3pt
{{#1}\atop{#2}}\kern-2pt \right)}
\def\smp#1#2#3{#1_{#2}+\dots+#1_{#3}} 

\def\lap#1{$\ell_{#1}$-$(\!${\it ap}$)$}\def\lpap/{\lap{p}}
\def\lmap#1{$\ell_{#1}$-$(\!${\it map}$)$}\def\lpmap/{\lmap{p}}
\def\uap/{\mbox{$(\!${\it uap}$)$}}\def\ubp/{\mbox{$(\!${\it ubp}$)$}}
\def\as/{\mbox{a.s.}}\def\sbd/{\mbox{$(\!${\it sbd\/}$)$}}
\def\ubs/{\mbox{$(\!${\it ubs}$)$}}\def\umbs/{\mbox{$(\!${\it umbs}$)$}}
\def\umap/{\mbox{$(\!${\it umap}$)$}}
\def\ap/{\mbox{$(\!${\it ap}$)$}}\def\map/{\mbox{$(\!${\it map}$)$}}
\def\fdd/{\mbox{$(${\it fdd\/}$)$}}\def\umfdd/{$(${\it umfdd\/}$)$}

\def\SC#1{\mathscr{C}_{#1}(\T)}

\def\SL#1#2{{\fam0 L}^{#1}_{#2}(\T)}

\def\SLP#1{\SL{p}{#1}}

\def\DX/{$X\in\{\SC{},\SLP{}:1\le p<\infty\}$}
\def\PT#1{\mbox{$\mathscr{P}_{#1}(\T)$}}
\def\PTE/{\PT{E}}

\def\UP/{\mbox{$(\mathscr{U})$}}

\def\EL#1{$\Lambda(#1)$}
\def\ER/{\mbox{Rosenthal}} 

\def\DE/{\mbox{$E\subseteq\Z$}}
\def\DEE/{\mbox{$E=\{n_k\}\subseteq\Z$}}
\def\DEEE/{\mbox{$E=\{n_k\}_{k\ge1}\subseteq\Z$}}
\def\DL/{\mbox{$\L/\subseteq\Z$}}
\def\DLL/{\mbox{$\L/=\{\lambda_k\}\subseteq\Z$}}
\def\DLLL/{\mbox{$\L/=\{\lambda_k\}_{k\ge1}\subseteq\Z$}}

\def\DI/{$I\subseteq\N\times\N$}

\def\ud/{equidistributed}  
\def\wud/{weakly equidistributed}  

\title{Two random constructions inside lacunary sets} 
\author{Stefan Neuwirth}
\date{}
\begin{document}

\def\set{E}
\def\element{n}
\def\SCL{\SC{\set}}
\def\SLL#1{\SL{#1}{\set}}
\def\SLPL{\SLL{p}}
\def\DL/{$\set\se\Z$}
\def\DLL/{$\set=\{\element_k\}\se\Z$}
\def\DLLL/{$\set=\{\element_k\}_{k\ge1}\se\Z$}

\maketitle

{\leftskip=1cm\rightskip=1cm\small

\begin{center}{\bf R\'esum\'e}\end{center}

Nous \'etudions le rapport entre la
croissance d'une suite d'entiers et les propri\'et\'es harmoniques et
fonctionnelles de la suite de caract\`eres associ\'ee. Nous montrons
en particulier que toute suite polynomiale, ainsi que la suite des
nombres premiers, contient un ensemble $\Lambda(p)$ pour tout $p$ qui
n'est pas de Rosenthal.

\begin{center}{\bf Abstract}\end{center}

We study the relationship between the growth rate
of an integer sequence and harmonic and functional properties of the
corresponding sequence of characters. We show in particular that every
polynomial sequence contains a set that is $\Lambda(p)$ for all $p$
but is not a Rosenthal set. This holds also for the sequence of
primes.

}

\section{An introduction in French}

\subsection{Position du probl\`eme} 

Nous \'etudions le rapport entre la croissance d'une suite
$\{n_k\}=E\se\Z$
et deux de ses propri\'et\'es harmoniques et fonctionnelles
\'eventuelles, \ie \\
\bloc 
toute fonction int\'egrable sur le tore \`a spectre dans $E$ est en fait
$p$-int\'egrable pour tout $p<\infty$: $E$ est un ensemble
$\Lambda(p)$ pour tout $p$;\\
\bloc 
toute fonction mesurable born\'ee sur le tore \`a spectre dans $E$ est continue
\`a un ensemble de mesure nulle pr\`es: $E$ est un ensemble de
Rosenthal.

Nous sommes en mesure de dresser le tableau suivant selon la
croissance\\
\bloc 
polynomiale: $n_k\preccurlyeq k^d$ pour un $d<\infty$, \\
\bloc 
surpolynomiale: $n_k\gg k^d$ pour tout $d\ge1$, \\
\bloc 
sous-exponentielle: $\log n_k\ll k$, \\
\bloc 
g\'eom\'etrique: $\liminf |n_{k+1}/n_k|>1$.

\begin{table}[h]
\label{intro:table}
  \hfil\vbox{\halign{\lower4pt\vbox to 14pt{}\hfil#\ \vrule\hbox to 1pt{}&\vrule\ 
  \hfil#\hfil\ &\vrule\ \hfil#\hfil\ &\vrule\ \hfil#\hfil\cr
  croissance&polynomiale&surpolynomiale et sous-exponentielle&g\'eom\'etrique\cr
\noalign{\hrule}\noalign{\vskip1pt\hrule}
  $E$ $\Lambda(p)$ $\forall p$& non& presque
  toujours&oui\cr\noalign{\hrule} $E$ Rosenthal&presque
  jamais\span\omit&oui\cr}}
\caption{Croissance et propri\'et\'es harmoniques ou fonctionnelles.}
\end{table}\medskip

Li montre qu'effectivement il existe un ensemble $\Lambda(p)$ pour
tout $p$ qui n'est pas de Rosenthal. Nous traitons les deux questions
suivantes. 

\begin{Q}
  Le sch\'ema ci-dessus reste-t-il valable si on consid\`ere \`a la
  place de l'ensemble des sous-ensembles $E$ de $\Z$ l'ensemble des
  sous-ensembles $E$ d'une suite \`a croissance polynomiale~?
\end{Q}

\begin{Q}
  \label{intro:random:q2}
  % ou aussi de la suite des nombres premiers. En particulier, il
  % existe \`a
  Si $E$ n'est pas un ensemble de Rosenthal, $E$ contient-il un
  ensemble \`a la fois $\Lambda(p)$ pour tout $p$ et non Rosenthal~?
\end{Q}

\subsection{Constructions al\'eatoires \`a l'int\'erieur de suites
  lacunaires} 

Nous fournissons une preuve nouvelle pour une
construction al\'eatoire d'ensembles lacunaires par Yitzhak Katznelson qui
appartient au folklore de l'analyse harmonique. Nous analysons et
g\'en\'eralisons aussi la construction al\'eatoire d'ensembles
\'equidistribu\'es par Jean Bourgain.

Cela nous permet d'\'etablir le tableau \ref{intro:table} qui classe
les propri\'et\'es de Rosenthal et $\Lambda(p)$ pour tout $p$ selon la
croissance du spectre. Nous montrons alors que la d\'emarche
probabiliste suivie par Katznelson et Bourgain pour construire ces
sous-ensembles de $\Z$ utilise seulement la croissance
``arithm\'etique'' et l'\'equidistribution de la suite des entiers
$\Z$. En fait, ces sous-ensembles peuvent \^etre construits \`a
l'int\'erieur de suites \'equidistribu\'ees \`a croissance
polynomiale. En particulier, le tableau \ref{intro:table} reste
valable pour l'ensemble des sous-ensembles $E$ d'une suite
polynomiale, ainsi que de la suite des nombres premiers.

Nous fournissons une r\'eponse partielle \`a la question
\ref{intro:random:q2}.
\begin{thmf}
Soit $P$ une suite polynomiale ou la suite des nombres
premiers. Alors il existe une sous-suite $E$ de $P$ qui est
$\Lambda(p)$ pour tout $p$ alors qu'elle ne forme pas un ensemble de
Rosenthal.
\end{thmf}

\section{Introduction}
The study of lacunary sets in Fourier analysis still suffers from a
severe lack of examples, in particular for the purpose of
distinguishing two properties. In order to bypass the individual
complexity of integer sets, one frequently resorts to random
constructions. In particular, Li \cite{li98} uses in his argumentation
a construction due to Katznelson \cite{ka73}  
%and a theorem of Bourgain \cite[Prop.\ 8.2(i)]{bo88} 
to discriminate the following two functional properties of certain
subsets \DL/: 

\vskipa
{\setbox0\hbox{\bloc\kern1ex}\parindent-\wd0\leftskip\wd0

\bloc\kern1ex
A Lebesgue integrable function on the circle with Fourier frequencies in 
$\set$ is in fact $p$-integrable for all $p<\infty$. This means that  
all spaces $\SLPL$ coincide for $p<\infty$, 
\ie $\set$ is a \EL{p}\index{Lambda(p) set@\EL{p} set} set for all
$p$ in Rudin's terminology. No sequence of polynomial growth has this
property \cite[Th.\ 3.5]{ru60}. By
Theorem \ref{lambda:katz-gen}, almost every sequence of a given
superpolynomial order of growth is \EL{p} for all $p$.

\vskipa\bloc\kern1ex
A bounded measurable function on the circle with Fourier frequencies in
$\set$ is in fact continuous up to a set of measure $0$. 
This means that $\SLL{\infty}$ and $\SCL$ coincide: $\set$ is a
\ER/\index{Rosenthal set} 
set. Every sequence of exponential growth is a Sidon set and therefore
has this property. By Bourgain's Theorem \ref{def:bourgain}, 
%\cite[Prop.\ 8.2(i)]{bo88}
almost every sequence
of a given subexponential order of growth fails the \ER/
property.\vskipb

}  

A Rosenthal set may contain arbitrarily
large intervals \cite{ro67} und thus fail the \EL{p} property. This
shows that these two properties cannot be characterized by some
order of growth, whereas the random method is so imprecise that it 
ignores a range of exceptional sets.  
On the other hand, Li shows that some set $E$ is \EL{p} for all $p$ and
fails the Rosenthal property: his 
construction witnesses for the quantitative overlap between
superpolynomial and subexponential order of growth. From a Banach
space point of view, Li's set $E$ is such that $\SCL$ contains $c_0$
while $\SLL{1}$ does not contain $\ell_1$.

We come back to Li \cite{li98} for two reasons: in the first place, 
we have been unable to locate a published proof of Katznelson's
statement. We
provide one for a 
stronger statement in Section \ref{lambda}. In the second place,  we
want to 
precise and supple the random construction in the following sense: can
one distinguish the \EL{p} property and the Rosenthal property among
subsets of a 
certain given set~? That sort of questions has been investigated
by Bourgain\index{Bourgain, Jean} in \cite{bo89a}. We give the
following answer (see Th.\ \ref{def:bourgain-gen}): 

\vskipa
{\bf Main Theorem } 
{\it Consider a polynomial\index{polynomial sequences} 
sequence of integers, or the sequence of
primes\index{prime numbers}. Then some subsequence of it is \EL{p} for
all $p$ and at the 
same time fails the Rosenthal property.}
\vskipb

This is a special case of the more general question: does every set that
fails the Rosenthal property contain a subset that is \EL{p} for all
$p$ and still fails the Rosenthal property~? We should emphasize at
this point that neither of these notions has an arithmetic
description. In fact, the family  of Rosenthal sets is coanalytic non
Borel \cite{go91} and any description would be at least as complex
as their definition. This is why we study instead the following two
properties for certain subsets \DL/.

\vskipa
{\setbox0\hbox{\bloc\kern1ex}\parindent-\wd0\leftskip\wd0

\bloc\kern1ex
Any integer $n$ has at most one representation as the sum of $s$
elements of $\set$. This implies that $\set$ is \EL{2s} by
\cite[Th.\ 4.5(b)]{ru60}. 

\vskipa\bloc\kern1ex
$\set$ is \ud/ in Hermann Weyl's\index{Weyl, Hermann} sense: save for
$t\equiv0$ mod $2\pi$, 
the 
successive means of $\{\e^{\ii\element t}\}_{\element\in\set}$ tend to
$0$, which is the mean of $\e^{\ii t}$ over $[0,2\pi\mathclose[$. This
implies that $\set$ is not a Rosenthal set by \cite[Lemma 4]{lu81}. \vskipb

}

Our random construction gives no hint for explicit procedures to build
such integer sets. The question whether some ``natural'' set of
integers is \EL{p} for all $p$ and fails the Rosenthal property
remains open. 

Let us describe the paper briefly. Section \ref{def} introduces the
inquired notions and gives a survey of former and new
results. As the right framework for this study appears to consist in
the sequences of polynomial growth, we give them a precise meaning in
Section \ref{poly}, and show that they are nicely distributed among
the intervals of the partition of \Z\ defined by
$\{\pm2^{k!}\}$. Section \ref{lambda} establishes an
optimal criterion for the generic subset of a set with polynomial
growth to be \EL{p} for all $p$. Section \ref{ergo} comes back to
Bourgain's proof in \cite[Prop.\ 8.2(i)]{bo88}: we simplify and
strengthen it in order to investigate the generic subset of an \ud/
set.

\vskipa{\bf Notation }
$\T=\{t\in\C:\,|t| = 1\}$ is the unit circle endowed with its Haar
measure $dm$ and \Z\ its dual group of
integers: for each $n\in\Z$, let $\e_n(t)=t^n$. The cardinal of \DLL/
is written \mes{\set}. We denote by $c_0(\T)$ the space of functions on \T\
which are arbitrarily small outside finite sets; such functions
necessarily have countable support.

For a space of integrable functions on \T\ and \DL/, $X_\set$
denotes the 
space of functions with Fourier spectrum in $\set$: 
$X_\set=\bigl\{f\in X: \widehat{f}(n)=\int\e_{-n}f\,dm=0$ if $n\notin\set\bigr\}$.

We shall stick to Hardy's notation: $u_n\preccurlyeq v_n$
(\vs $u_n\ll 
v_n$) if $u_n/v_n$ is bounded (\vs vanishes) at infinity.

\vskipa{\bf Acknowledgment } 
I would like to thank Daniel Li for several helpful discussions.
\section{\texorpdfstring{Equidistributed and \EL{p} sets}{Equidistributed and Lambda(p) sets}}\label{def}
\begin{dfn}
Let \DLLL/ ordered by increasing absolute value $|\element_k|$. 

$(i)$
\cite[Def.\ 1.5]{ru60}
Let $p>0$. $\set$ is a \EL{p} set if, for some ---~or equivalently for any~---
$0<r<p$, $\SLPL$ and $\SLL{r}$ coincide:
$$
\exists C_r\quad\forall f\in\SLPL\quad\|f\|_r\le\|f\|_p\le C_r\|f\|_r.
$$

$(ii)$ 
\cite[\S 7]{we16}
$\set$ is \ud/\index{equidistributed set} if for each
$t\in\T\setminus\{1\}$ the 
successive means 
\begin{equation}\label{def:def}
f_k(t)=\frac1k\sum_{j=1}^k\e_{\element_j}(t)\ \tol_{k\to\infty}\ 0.
%\hbox{\fam0\bf 1}_{\{1\}}(t)=\left\{\begin{array}{l}
%      \mbox{$1$ if $t=1$}\\
%      \mbox{$0$ if not.}
%\end{array}\right.
%
\end{equation}
Thus $E$ is equidistributed if and only if the sequence of characters
in $E$ converges to $\hbox{\fam0\bf 1}_{\{1\}}$ for the Ces\`aro
summing method. 
If $f_k$ tends pointwise to $f\in c_0(\T)$, then $\set$ is
\wud/\index{weakly equidistributed set}.
\end{dfn}
If $\set$ is \wud/, then $f$ defines an element of
$\SCL^{\perp\perp}$.  By Lust-Piquard's \cite[Lemma 4]{lu81}, $\SCL$
then contains a copy of $c_0$ and $\set$ cannot be a Rosenthal set. 
%See also \cite{lu89} for further background.

For example, \Z\ and \N\ are \ud/.  Arithmetic\index{arithmetic
  sequences} sequences are \wud/: there is a finite set on which
$f_k\nrightarrow0$. Polynomial sequences\index{polynomial sequences}
of integers (\cite[Th.\ 9]{we16} and \cite[Lemma 2.4]{va81}, 
see \cite[Ex.\ 2]{lu89}) 
and the sequence of prime\index{prime numbers} numbers 
(Vinogradov's\index{Vinogradov} theorem \cite{el75}, see \cite[Ex.\ 1]{lu89}) 
are \wud/: $f_k(t)$
may not converge to $0$ for rational $t$. There are nevertheless
sequences of bounded pace that are not \wud/ \cite[Th.\ 11]{et57}.
Sidon sets are \EL{p} for all $p$ \cite[Th.\ 3.1]{ru60}, but not \wud/
since they are Rosenthal sets. \vskipb

\exa
Consider the geometric\index{geometric sequences}
sequence $\set=\{3^k\}_{k\ge1}$ and the corresponding sequence of
successive means $f_k$. By
\cite[Th.\ 14]{et57}, the $f_k$ do not converge to 0 on a null set of
Hausdorff dimension $1$. Consider 
$$
f_k^j=k^{-j}\hskip-1.5em\sum_{1\le\lst kj\le k}\hskip-1.5em
\e_{3^{k_1}+\dots+3^{k_j}}=
k^{-j}\Bigl(\hskip-1ex\mathop{j!\sum}\limits_{1\le k_1<\dots<k_j\le k}+
\sum_
{\substack{1\le\lst kj\le k\\\text{not all distinct}}}
\hskip-1ex\Bigr)
\e_{3^{k_1}+\dots+3^{k_j}}.
$$
Let $j\ge1$. Put 
$\set^{(j)}=\{3^{k_1}+\dots+3^{k_j}:\,1\le k_1<\dots<k_j\}$
and let $f_k^{(j)}$ be the corresponding successive means
\eqref{def:def}. Then
\begin{eqnarray*}
\|f_k^j-f_{\binom kj}^{(j)}\|_\infty&\le&
\Bigl(\binom kj^{-1}-\frac{j!}{k^j}\Bigr)
\binom kj+\frac1{k^j}
\Bigl(k^j-\frac{k!}{(k-j)!}\Bigr)\\
&=&2\Bigl(1-\frac{k!}{k^j(k-j)!}\Bigr)
\tol_k0.
\end{eqnarray*}
Therefore the $f_k^{(j)}$ do not converge to $0$ outside a countable set,
and $\set^{(j)}$, which is \EL{p} for all $p$ \cite[Th.\ IV.3]{me68}
and not Sidon, is not \wud/. \vskipb

However, as Li\index{Li, Daniel} notes, these two classes meet.
\begin{thm}[\cite{li98}]\label{def:li}
There is an \ud/ sequence that is \EL{p} for all $p$.
\end{thm}
{\it Sketch of proof.}\/ 
Li uses the following random construction, discovered by
Erd\H{o}s\index{Erd\H{o}s, Paul} \cite{er56,er60} and introduced to
harmonic analysis by 
Katznelson\index{Katznelson, Yitzhak} and 
Malliavin\index{Malliavin, Paul} \cite{km65,km66}.
\begin{construction}\label{def:constr}
\index{random construction}
Let \DL/ and consider independent $\{0,1\}$-valued selectors
$\xi_\element$ of mean $\delta_\element$ $(\element\in\set)$, \ie
$\P{\xi_\element=1}=\delta_\element$. Then the random set $\set'$ is
defined by
$$
\set'=\{\element\in\set:\ \xi_\element=1\}.
$$
\end{construction}
The first ingredient of the proof is Bourgain's\index{Bourgain, Jean}
following  
\begin{thm}[\protect{\cite[Prop.\ 8.2(i)]{bo88}}]\label{def:bourgain}
%\index{equidistributed set}
Let $\set=\N$ in Construction \ref{def:constr}. If
$\delta_\element$ decreases with $\element$ while
$\delta_\element\gg\element^{-1}$, then $\set'$ is almost surely \ud/.
%skip
\end{thm}
\rem
In this sense, almost every sequence of a subexponential 
growth\index{subexponential growth} 
given
by $\{\delta_\element\}$ is equidistributed: indeed, for almost every
$E'\se\N$, 
$$
\mes{\set'\cap[0,\element]}
\sim\smp\delta0\element\gg\log\element
$$
by the Law of Large Numbers. Note however that the set $\set^{(j)}$
defined in Example 2.2 has subexponential growth:
$\mes{\set^{(j)}\cap[-\element,\element]}\succcurlyeq(\log\element)^j$, and
is not \ud/.\vskipb

The second ingredient is a result announced without proof by
Katznelson.
\begin{prp}[\protect{\cite[\S 2]{ka73}}]\label{def:katz}
Put $I_k=\mathopen]p_{k-1},p_k ]$ with $p_k >p_{k-1}^2$ $(k\ge1)$. Let
$\set=\N$ in Construction \ref{def:constr}. There is a choice of
$(\ell_k)$ with $\ell_k\gg\log p_k$ such that for
$\delta_\element=\ell_k/\lmes I_k\rmes $ $(\element\in I_k)$, $\set'$
is \EL{p} for all $p$ almost surely.
\end{prp}
Li suggests to apply the content of Proposition \ref{def:katz} with
$p_k=2^k$ and $\ell_k=k$: then $\delta_\element\gg\element^{-1}$ and
Theorem \ref{def:li} derives from Theorem \ref{def:bourgain}.
\eck\vskipb

We shall generalize Katznelson's and Li's results with a new proof that
permits to construct $\set'$ inside of sets $\set$ with polynomial growth
(see Def.\ \ref{poly:def}) and yields an optimal criterion on
$\ell_k$. We shall subsequently generalize Bourgains's Theorem
\ref{def:bourgain} 
to obtain the Main Theorem via
\begin{thm}\label{def:bourgain-gen}
Let $\set$ be \ud/ \(\vs weakly\) and with polynomial
growth\index{polynomial growth}.
Then there is a subset $\set'\subseteq\set$ \ud/ \(\vs weakly\)
and at the same time \EL{p} for all $p$.
\end{thm}
A precise and quantitative statement of this is Theorem
\ref{ergo:COR}.
\section{Sets with polynomial growth}\label{poly}
We start with the definition and first property of such sets.
\begin{dfn}\label{poly:def}
Let \DLLL/ be an infinite set ordered by increasing absolute value and
$\set[t]=\lmes\set\cap[-t,t]\rmes $ its distribution 
function\index{distribution function}. 

$(i)$
$\set$ has polynomial growth\index{polynomial growth} 
if $\element_k\preccurlyeq k^d$ for
some $1 \le d < \infty$. This amounts to $\set[t]\succcurlyeq t^\eps$
for $\eps=d^{-1}$.

$(ii)$
$\set$ has regular polynomial growth\index{regular polynomial growth} 
if there is a $c>1$ such that 
$|\element_{\lceil ck\rceil}| \le 2|\element_k|$ for large $k$.
This amounts to $\set[2t] \ge c\set[t]$ for large $t$.
\end{dfn}
\dem
$(i)$
If $|\element_k| \le  C k^d$  for large $k$ and $C k^d \le  t < C (k+1)^d$,
then $\set[t] \ge k > (t/C)^\eps-1$. Conversely, if $\set[t] \ge c
t^\eps$ for large $t$ and $c (t-1)^\eps < k \le c t^\eps$, then
$|\element_k| \le  t < (k/c)^d+1$. 

$(ii)$ If $|\element_{\lceil ck\rceil}| \le  2|\element_k|$ for large $k$
and $k$ is maximal with $|\element_k|\le t$, then $\set[2t] \ge
\set[2|\element_k|] \ge c k = c \set[t]$. Conversely, if $\set[2t]
\ge c\set[t]$ for large $t$, then $\set[|\element_k|] \in \{k,k+1\}$ and 
$\set[2|\element_k|]\ge c k$. Thus $|\element_{\lceil ck\rceil}| \le
2|\element_k|$. 
\eck\vskipb

In particular, polynomial sequences\index{polynomial sequences} 
have regular polynomial growth. By
the Prime Number Theorem, the sequence of primes\index{prime numbers} 
also has.
Property $(ii)$ implies property $(i)$: if $\set[2t] \ge c\set[t]$ for
large $t$, then $\set[t]\succcurlyeq t^{\log_2c}$. The converse
however is false as shows
$F=\bigcup\, \mathopen]2^{2^{2k}},2^{2^{2k}+1}]$, for which $F[t]
\succcurlyeq t^{1/4}$ while $F[2t]=F[t]$ infinitely
often. 

Let us relate Definition \ref{poly:def} with certain partitions
of \Z. Regular growth means in fact that $\set$ is regularly distributed 
on the annular dyadic partition\index{dyadic partition} of \Z
\begin{equation}\label{poly:lpd}
\mathscr{P}=\bigl\{[-p_0,p_0],I_k=[-p_k,-p_{k-1}\mathclose[\cup
\mathopen]p_{k-1},p_k]\bigr\}_{k\ge1}\mbox{ where }p_k=2^k
\end{equation} 
and $F$ shows that there are sets with polynomial growth which are not
regularly distributed on the partition defined by
$p_k=2^{2^k}$. However, the intervals of the gross 
partition\index{gross partition} 
\begin{equation}\label{poly:lpg}
\mathscr{P}=\bigl\{[-p_0,p_0],I_k=[-p_k,-p_{k-1}\mathclose[
\cup\mathopen]p_{k-1},p_k]\bigr\}
\mbox{ where }\log p_k \gg\log p_{k-1}
\end{equation} 
grow with a speed that forces regularity. Put $p_k=2^{k!}$ for a
simple explicit example. We have precisely 
\begin{prp}\label{poly:regul}
Let \DL/, $\mathscr{P}=\{I_k\}$ a partition of \Z\ and $\set_k=\set\cap
I_k$. Then $\log\mes{\set_k}\succcurlyeq \log \mes{I_k}$ in the two
following cases:

$(i)$
$\set$ has regular polynomial growth and $\mathscr{P}$ is partition
\eqref{poly:lpd};

$(ii)$
$\set$ has polynomial growth and $\mathscr{P}$ is partition
\eqref{poly:lpg}.
\end{prp}
\dem
$(i)$ Choose $K$ and $c>1$ such that $\set[2^k] \ge c\set[2^{k-1}]$
for $k\ge K$. Then $\set[2^k] \succcurlyeq c^k$.
%if \set is non-empty. 
Thus
$$
\mes{\set_k} = \set[2^k] - \set[2^{k-1}] \ge  
(1-c^{-1})\set[2^k] \succcurlyeq c^k=2^{k\log_2c}.
$$
$(ii)$
In this case $p_k ^\eps \gg p_{k-1}$ for each positive $\eps$. Now there is
$\eps>0$ such that 
$$
\mes{\set_k} = \set[p_k] - \set[p_{k-1}] \succcurlyeq p_k ^\eps \succcurlyeq
\mes{I_k}^\eps.\eqno{\eck}
$$
\section{\texorpdfstring{Sets that are \EL{p} for all $p$}{Sets that are Lambda(p) for all p}}\label{lambda}
In this section, we establish an improvement (Th.\
\ref{lambda:katz-gen}) of Katznelson's statement \cite[\S2]{ka73}. We
  first recall several known definitions and results. 

\EL{p} sets have a practical description in terms of
unconditionality. We shall also use a combinatorial property that is
more elementary than \cite[1.6(b)]{ru60}: to this end, write
$\Zeta_s^m$ for the following set of arithmetic relations.
$$
\Zeta_s^m=\bigl\{\zeta\in{\Z^*}^m:\,\sm\zeta m=0\mbox{ and }
|\zeta_1|+\dots+|\zeta_m|\le 2s\bigr\}. 
$$
Note that $\Zeta_s^1$ and $\Zeta_s^m$ $(m>2s)$ are empty, and that
every $\zeta\in\Zeta_s^2$  is of form $\zeta_1\cdot(1,-1)$: this is
the identity relation.
\begin{dfn}
Let $1\le p<\infty$, $s\ge1$ integer and \DL/.

$(i)$ \cite{ka48}
$\set$ is an unconditional basic sequence in $\SLP{}$ if 
$$
\sup_{\pm}\biggl\|\sum_{n\in\set}\pm a_\element\e_\element\biggr\|_p\le
C\biggl\|\sum_{\element\in \set}a_\element\e_\element\biggr\|_p.
$$
for some $C$. If $C=1$ works, $\set$ is a $1$-unconditional basic
sequence in $\SLP{}$.

$(ii)$
%[Chapter \ref{chap1}, \S\ref{ss:isom}] 
$\set$ is $s$-independent if $\sum_1^m\zeta_iq_i\ne0$ for all $3\le m\le 2s$, 
$\zeta\in\Zeta_s^m$ and distinct $\lst qm\in\set$.
\end{dfn}
\begin{prp}
%[cf.\ \protect{\cite[Prop.\ 2.5]{ne98}}]
\label{lambda:uncond}
Let $1\le p\ne2<\infty$, $s\ge1$ integer and \DL/.

$(i)$ 
\cite[proof of Th.\ 3.1]{ru60}
$\set$ is a \EL{\max(p,2)} set if and only if $\set$ is an
unconditional basic sequence in $\SLP{}$.

$(ii)$
%[Chapter \ref{chap1}, Prop.\ \ref{mub:isom}]
\cite[Prop.\ 2.5, Rem.\ (1)]{ne98}
$\set$ is a $1$-unconditional basic sequence in $\SL{2s}{}$ if and only
if $\set$ is $s$-independent.
\end{prp}
We need to introduce a second classical notion of unconditionality
that rests on the Little\-wood--Paley theory.
\begin{dfn}[\cite{hk89}]
Let $\mathscr{P }=\{I_k\}$ be a partition of \Z\ in finite
intervals. $\mathscr{P}$ is a Littlewood--Paley
partition\index{Littlewood--Paley partition} if for each
$1<p<\infty$ there is a constant $C_p$ such that 
\begin{equation}\label{lambda:paley}
\forall f\in\SLP{}\quad
\sup_{\pm}\Bigl\|\sum\pm f_k\Bigr\|_p\le C_p\|f\|_p\qquad
\mbox{with }\widehat{f_k}=
\left\{\begin{array}{l}\widehat{f}\mbox{ on}\\
0\mbox{ off}\end{array}\right. I_k. 
\end{equation}
\end{dfn}
By Khinchin's inequality, this means exactly that 
$$
\forall f\in\SLP{}\quad
\|f\|_p\approx\Bigl\|\Bigl(\sum|f_k|^2\Bigr)^{1/2}\Bigr\|_p.
$$
In particular, the dyadic partition \eqref{poly:lpd} and the
gross partition \eqref{poly:lpg} are Littlewood--Paley
\cite{lp31}. By Proposition \ref{lambda:uncond} and
\eqref{lambda:paley}, we obtain 
\begin{prp}\label{lambda:union}
Let $\{I_k\}$ be a Littlewood--Paley partition and $\set_k\subseteq
I_k$. If $\set_k$ is $s$-independent for each $k$, then
$\set=\bigcup\set_k$ is an unconditional basic sequence in $\SL{2s}{}$ and
thus a \EL{2s} set.
\end{prp}
We generalize now Katznelson's Proposition \ref{def:katz}.
\begin{lem}\label{lambda:estim}
Let $s\ge2$ integer, \DL/ finite and $0\le\ell\le\mes{\set}$. Put
$\delta_n=\ell/\mes{E}$ in
Construction \ref{def:constr}, so that all selectors $\xi_\element$
have same distribution. Then there is a constant $C(s)$ that depends
only on $s$ such that
$$
\P{\set'\mbox{ is $s$-dependent\/}} \le 
C(s)\frac{\ell^{2s}}{\lmes \set\rmes }.
$$
\end{lem}
\dem
We wish to compute the probability that there are $3\le m\le 2s$,
$\zeta\in\Zeta_s^m$ and distinct $\lst qm\in\set'$ with
$\sum\zeta_iq_i=0$. As the number $C(s)$ of arithmetic relations 
$\zeta\in\Zeta_s^m$ ($3\le m\le 2s$) is finite and depends on $s$
only, it suffices to compute, for fixed $m$ and $\zeta\in\Zeta_s^m$
\begin{eqnarray*}
\lefteqn{\P{\exists\lst qm\in\set'\mbox{ distinct}:\,
\sum\zeta_iq_i=0}}\\
&=&\P{\exists\lst q{m-1}\in\set'\mbox{ distinct}:\,
-\zeta_m^{-1}\sum_{i=1}^{m-1}\zeta_iq_i\in
\set'\setminus\{\lst q{m-1}\}}\\
&=&\P{\bigcup_
{\substack{\lst q{m-1}\\\in\set'\,\text{distinct}}}\hskip-1.5em
\Bigl[-\zeta_m^{-1}\sum_{i=1}^{m-1}\zeta_iq_i\in
\set'\setminus\{\lst q{m-1}\}\Bigr]}\\
&=&\P{\bigcup_
{\substack{\lst q{m-1}\\
\in\set\,\text{distinct}}}\hskip-1.5em
\Bigl[q_m=-\zeta_m^{-1}\sum_{i=1}^{m-1}\zeta_iq_i\in
\set\setminus\{q_i\}_{i=1}^{m-1}\ \&\ 
\xi_{q_1}=\dots=\xi_{q_m}=1\Bigr]}.
\end{eqnarray*}
The union in the line above runs over 
$$
\frac{\lmes \set\rmes !}{(\lmes \set\rmes -m+1)!}\le\lmes \set\rmes ^{m-1}
$$
$(m-1)$-tuples. Further, the event in the inner brackets implies that
$m$ out of 
$\mes{\set}$ selectors $\xi_\element$ have value $1$: its probability is
bounded by $(\ell/\mes{\set})^m$. Thus
$$
\P{\set'\mbox{ is $s$-dependent}}\le
C(s)\max_{3\le m\le 2s}\lmes \set\rmes ^{m-1}\frac{\ell^m}{\lmes
\set \rmes ^m}\le
C(s)\frac{\ell^{2s}}{\lmes \set\rmes }.\eqno{\eck}
$$
The random method we shall use is the following random construction. 
\begin{construction}\label{lambda:constr}\index{random construction}
Let \DL/. Let $\{I_k\}$ be a Littlewood--Paley partition and
$\set_k=\set\cap I_k$. Let $(\ell_k)_{k\ge1}$ with $0\le\ell_k\le\lmes
\set_k\rmes $ and put 
$$
\P{\xi_\element=1}=\delta_\element=\ell_k/\mes{\set_k}\quad
(\element\in\set_k)
$$
in Construction \ref{def:constr}. Put $\set'_k=\set'\cap I_k$.
\end{construction}
\begin{thm}\label{lambda:katz-gen}
Let \DL/ have polynomial \(\vs regular\) growth and $\{I_k\}$ be the
gross \eqref{poly:lpg} \(\vs dyadic \eqref{poly:lpd}\)
Littlewood--Paley partition. Do Construction \ref{lambda:constr}. The
following assertions are equivalent.

$(i)$
$\log \ell_k \ll \log \lmes I_k \rmes $, \ie $\log \ell_k \ll \log
p_k $ \(\vs $\log \ell_k \ll k$\);

$(ii)$
$\set'$ is almost surely a \EL{p}\index{Lambda(p) set@\EL{p} set!for all $p$} 
set for all $p$. 
\end{thm}
\dem Note that by Proposition \ref{poly:regul}, there is a positive
$\alpha$ such that $\mes{\set_k}>\mes{I_k}^{\alpha}$ for large $k$.

$(i)\imp(ii)$ Let $s\ge2$ be an arbitrary integer. By Proposition 
\ref{lambda:estim}, 
$$
\sum_{k=1}^\infty\P{\set'_k\mbox{ is $s$-dependent}}\le
C(s)\sum_{k=1}^\infty\frac{\ell_k^{2s}}{\lmes \set_k\rmes }.
$$
For each $\eta>0$, $\ell_k\le\mes{I_k}^\eta$ for large $k$. Choose
$\eta<\alpha/2s$. Then  $\ell_k^{2s}/\mes{\set_k}\le
\mes{I_k}^{2s\eta-\alpha}$  for large
$k$, and the series above converges since $\mes{I_k}\succcurlyeq
2^k$. By the Borel--Cantelli lemma, $\set'_k$ is almost surely
$s$-independent for large $k$. By Proposition \ref{lambda:union}, $\set'$ is
almost surely the union of a finite set and a \EL{2s} set. By
\cite[Th.\ 4.4(a)]{ru60}, $\set'$ itself is almost surely a \EL{2s} set.

$(ii)\imp(i)$ 
If $\set'$ is a \EL{2s} set, then by \cite[Th.\ 3.5]{ru60} or simply by
\cite[(1.12)]{bo89a}, there is a constant $C_s$ such that
$\mes{\set'_k}<C_s\mes{I_k}^{1/s}$. As $\mes{\set'_k}\sim\ell_k$ almost
surely by the Law of Large Numbers (cf.\ the following Lemma
\ref{ergo:bern}), one has $\log \ell_k \ll \log \mes{I_k}$.
\eck\vskipb

\rem
As one may easily construct sets that grow as slowly as one wishes and
nevertheless contain arbitrarily large intervals 
(see also \cite[Th.\ 3.8]{ru60}
for an optimal statement), one cannot remove
the adverb ``almost surely'' in Theorem
\ref{lambda:katz-gen}$(ii)$.\vskipb 

\rem
The right formulation of Katznelson's Proposition \ref{def:katz} 
thus turns out to be the following. Let $\set=\N$ and
$I_k=\mathopen]p_{k-1},p_k ]$ with $p_k >cp_{k-1}$ for some $c>1$ in
Construction \ref{lambda:constr}. Then $\set'$ is almost surely a \EL{p}
set for all $p$ if and only if $\log \ell_k \ll \log p_k$. \vskipb 

\rem
Theorem \ref{lambda:katz-gen} shows that there are sets that are
\EL{p} for all $p$ of any  given
superpolynomial\index{superpolynomial growth} order of growth. This is
optimal since sets with 
distribution $\set[t]\succcurlyeq t^\eps$ fail the \EL{p} property for
$p>2/\eps$ by \cite[Th.\ 3.5]{ru60}. Such sets may also be constructed
inductively by combinatorial means: see 
%Section \ref{chap1}.\ref{sect:comb} and 
\cite[\S II, (3.52)]{hr83}. 
%\cite[\S10]{ne98}.
%
%
%
%
%
\section{Equidistributed sets}\label{ergo}
In this section, we shall finally state and prove our principal
result. To this end, we shall first generalize Bourgain's Theorem
\ref{def:bourgain} in order to get 
Theorem \ref{ergo:bourgain-gen}. 

The following lemma is Bernstein's
distribution inequality\index{Bernstein, Serge} \cite{be64} and dates
back to 1924. 
\begin{lem}\label{ergo:bern}
Let $\lst Xn$ be complex independent random variables such that
\begin{equation}\label{ergo:bern:hyp}
|X_i|\le1\quad\mbox{and}\quad\E X_i=0\quad\mbox{and}\quad
\E|X_1|^2+\dots+\E|X_n|^2\le\sigma.
\end{equation}
Then, for all positive $a$,
\begin{equation}\label{ergo:bern:conc}
\P{|X_1+\dots+X_n| \ge a}<4\exp(-a^2/4(\sigma+a)).
\end{equation}
\end{lem}
\dem
%This is Bennett's variant of Bernstein's 
%probabilistic inequality \cite{be64}. 
Consider first the case of real random variables. By
\cite[(8b)]{be62}, 
$$
\P{X_1+\dots+X_n \ge a}<\exp(a-(\sigma+a)\log(1+a/\sigma));
$$
as $\log(1-u)\le-u-u^2/2$ for $0 \le u < 1$,
$$
\P{X_1+\dots+X_n \ge a}<\exp(-a^2/2(\sigma+a)).
$$
One gets \eqref{ergo:bern:conc} since for complex $z$
$$
|z| \ge a\quad\Longrightarrow\quad
\max(\Re z,-\Re z,\Im z,-\Im z) \ge a/\sqrt{2}.\eqno{\eck}
$$
The next lemma corresponds to \cite[Lemma 8.8]{bo88} and is the
crucial step in the estimation of the successive means of 
$\{\e^{\ii\element t}\}_{\element\in\set'}$. Note that its hypothesis
is not on the individual $\delta_\element$, but on their successive
sums $\sigma_k$: this is needed in order to cope with the irregularity
of $\set$.
\begin{lem}\label{ergo:bourgain-lem}
Let \DLL/ be ordered by increasing absolute value. Do Construction
\ref{def:constr} and put
$\sigma_k=\delta_{\element_1}+\dots+\delta_{\element_k}$. If $\sigma_k
\gg \log |\element_k|$, then almost surely
\begin{equation}\label{ergo:bl:conc}
\psi(k)=\biggl\|
\frac1{\lmes \set'\cap\{\lst\element k\}\rmes }
\mathop{\sum\ \e_\element}\limits_{\set'\cap\{\lst\element k\}}-
\frac 1{\sigma_k}\sum_{j=1}^k\delta_{\element_j}\e_{\element_j}
\biggr\|_\infty\tol_{k\to\infty}0.
\end{equation}
\end{lem}
\dem
Note that 
$$
\mathop{\sum\ \e_\element}\limits_{\set'\cap\{\lst\element k\}}= 
\sum_{j=1}^k\xi_{\element_j}\e_{\element_j}\quad,\quad
\mes{\set'\cap\{\lst\element k\}}
=\sum_{j=1}^k\xi_{\element_j}.
$$ 
Center the $\xi_\element$ by putting 
$f=\sum_{j=1}^k(\xi_{\element_j}-\delta_{\element_j})\e_{\element_j}$. Then
\begin{eqnarray*}
\psi(k)&\le&\Bigl\|
\Bigl(\mes{\set'\cap\{\lst\element k\}}^{-1}-\sigma_k^{-1}\Bigr)
\sum_{j=1}^k\xi_{\element_j}\e_{\element_j}\Bigr\|_\infty+
\|\sigma_k^{-1}f\|_\infty\\
&\le&
\sigma_k^{-1}\,\biggl|\frac {\delta_{\element_1}+\dots+\delta_{\element_k}}
{\xi_{\element_1}+\dots+\xi_{\element_k}} - 1\biggr|\ 
\sum_{j=1}^k\xi_{\element_j}+
\sigma_k^{-1}\|f\|_\infty
\le2\sigma_k^{-1}\|f\|_\infty.
\end{eqnarray*}
Let $R=\{t\in\T:\,t^{4|\element_k|}=1\}$ and $u\in\T$ such that
$|f(u)|=\|f\|_\infty$. Let $t\in R$ be at
minimal distance of $u$: then $|t-u|\le\pi/4|\element_k|$. By
Bernstein's theorem,
$$
\|f\|_\infty-|f(t)| \le |f(u)-f(t)| \le |t-u| \, \|f'\|_\infty \le
\frac45 \|f\|_\infty;
$$
$$
\|f\|_\infty\le5\sup_{t\in R}|f(t)|.
$$
(For an optimal bound, cf.\ \cite[\S I.4, Lemma 8]{me72}.) 
%As $\mes{R}=4|\element_k|$ and 
For each $t\in R$, the random variables
$X_j=(\xi_{\element_j}-\delta_{\element_j})\e_{\element_j}(t)$ satisfy
\eqref{ergo:bern:hyp}, so that
$$
\P{|f(t)| \ge a}<4\exp(-a^2/4(\sigma_k+a)).
$$ 
As $\mes{R\mkern1mu}=4|\element_k|$, 
$$
\P{\|f\|_\infty \ge 5a}\le\P{\sup_{t\in
  R}|f(t)| \ge a}<4|\element_k|\cdot4\exp(-a^2/4(\sigma_k+a)).
$$
Put $a_k=(12\sigma_k\log |\element_k|)^{1/2}$. Then
$a_k\ll\sigma_k$: therefore
%: let $N$ such that $a_k\le\sigma_k$ pour $k\ge N$. Donc
$$
\P{\|f\|_\infty \ge 5a_k} \preccurlyeq |\element_k|^{-2}
$$
and by the Borel--Cantelli lemma, 
$$
\sigma_k^{-1}\|f\|_\infty\preccurlyeq a_k/\sigma_k\tol0\hbox{ almost surely}.
\eqno{\ecks}
$$

\rem
The hypothesis in Lemma \ref{ergo:bourgain-lem} contains implicitly a
restriction on the lacunarity of $\set$. If $\sigma_k \gg \log
|\element_k|$, then necessarily $\log|\element_k| \ll k$ and $\set[t]
\gg \log t$. In particular, $\set$ cannot be a Sidon set by
\cite[Cor.\ of Th.\ 3.6]{ru60}. 

\vskipa
We now state and prove the equidistributed counterpart of Theorem
\ref{lambda:katz-gen}.
\begin{thm}\label{ergo:bourgain-gen}
Let \DLL/ be equidistributed\index{equidistributed set} 
\(\vs weakly\index{weakly equidistributed set}\), and ordered by
increasing absolute value. Do Construction \ref{def:constr} and
suppose that $\delta_{\element_j}$ decreases with $j$. Put
$\sigma_k=\delta_{\element_1}+\dots+\delta_{\element_k}$. If $\sigma_k
\gg \log|\element_k|$, then $\set'$ is almost surely equidistributed
\(\vs weakly\). This is in particular the case if

$(a)$
$\ds\delta_{\element_k} \gg (|\element_k|-|\element_{k-1}|)/
{|\element_{k-1}|}$; 

$(b)$ 
$\set$ has polynomial growth and $\delta_{\element_k}\gg k^{-1}$.
\end{thm}
\dem
Lemma \ref{ergo:bourgain-lem} shows that almost surely
\eqref{ergo:bl:conc} holds. It remains to show that 
$$
\lim \frac 1 {\sigma_k}\sum_{j=1}^k\delta_{\element_j}\e_{\element_j} = \lim
\frac 1{k} \sum_{j=1}^k\e_{\element_j}, 
$$
\ie that the matrix summing method\index{matrix summing method} 
$(a_{k,j})$ given by 
$$
a_{k,j}=\begin{cases}\delta_{\element_j}/\sigma_k&\text{if $j\le k$}\\0&\text{if not}
\end{cases}
$$
is regular and stronger than the Ces\`aro 
method\index{Ces\`aro summing method} $C_1$ by arithmetic
means. This is the case because $a_{k,j}\ge0$,
$\sum_ja_{k,j}=1$ and (cf.\ \cite[\S52, Th.\ I]{ze58}) 
$$
\forall k\quad\sum_j j|a_{k,j}-a_{k,j+1}|=
\sum_j j(a_{k,j}-a_{k,j+1})=1<\infty
$$
since $a_{k,j}$ decreases with $j$ for each $k$.

$(a)$ In this case $\delta_{\element_k} \gg
\log|\element_k|-\log|\element_{k-1}|$ and thus $\sigma_k \gg
\log|\element_k|$. 

$(b)$
In this case, $\sigma_k \gg \log k \succcurlyeq \log
|\element_k|$.
\eck\vskipb 

In conclusion, we obtain, by combining Theorems \ref{lambda:katz-gen} and
\ref{ergo:bourgain-gen}, our principal result.
\begin{thm}\label{ergo:COR}
Let \DL/ be \ud/ \(\vs weakly\) and do Construction
\ref{lambda:constr}. Then
$\set'$ is almost surely \EL{p}\index{Lambda(p) set@\EL{p} set!for all $p$}  for all
$p$ and at the same time \ud/\index{equidistributed set}
\(\vs weakly\index{weakly equidistributed set}\) in the two following cases:

$(i)$
$\set$ is a set of regular polynomial growth, $\{I_j\}$ is the dyadic
Littlewood--Paley partition \eqref{poly:lpd}, $1 \ll \log\ell_j \ll j$ and
$\ell_j/\mes{\set_j}$ decreases eventually. 

$(ii)$
$\set$ is a set of polynomial growth, $\{I_j\}$ is the gross
Littlewood--Paley partition \eqref{poly:lpg}, $\ell_j/\mes{\set_j}$
decreases eventually 
and $\ell_j 
\gg \log p_{j+1}$ while 
$\log\ell_j \ll \log p_j $. This is the case if we put 
$p_j=2^{j!}$ and $\ell_j=\min((j+2)!,\mes{\set_j})$.
\end{thm}
\dem
In each case $\log\ell_j \ll \log\lmes
I_j\rmes$. Let us show that the hypothesis of Theorem
\ref{ergo:bourgain-gen} is verified. If
$\element_k\in\set_j\subseteq I_j$, then $|\element_k|\le p_j $ and 
$$
\sigma_k\ge\sum_{i=1}^{j-1}\sum_{\element\in\set_i}\delta_\element=
\sm\ell{j-1}
$$
and in each case $\ell_{j-1} \gg \log p_j -\log p_{j-1}$.

Let us make sure in $(ii)$ that our choice for $p_j$ and $\ell_j$
is accurate. Indeed, there is an $\eps>0$ such that
$\mes{\set_j}\succcurlyeq2^{\eps j!}$. Thus $(j+2)!\ll\mes{\set_j}$ and
$\ell_j=(j+2)!$ for large $j$. Note further that $(j+2)!\gg(j+1)!$ while
$\log(j+2)!\ll j!$. Finally
$$
\frac {\ell_{j+1}} {\mes{\set_{j+1}}} \preccurlyeq \frac {(j+3)!}
{2^{\eps(j+1)!}}  \preccurlyeq \frac {j\ell_j}{2^{\eps(j+1)!}}
\ll \frac {\ell_j} {2^{j!}} \preccurlyeq \frac {\ell_j} {\mes{\set_j}},
$$
so that $\ell_j/\mes{\set_j}$ decreases eventually.\eck

\def\cprime{$'$}\def\bibmath{disponible \`a la biblioth\`eque de
  math\'ematiques}\def\bibirem{disponible \`a la biblioth\`eque de
  l'IREM}\def\BUL{disponible \`a la BU Lettres}\def\BUS{disponible \`a la BU
  Sciences}\ifx\iflanguage\undefined\def\iflanguage#1#2#3{#3}\fi

\end{document}